\input amstex.tex
\input amsppt.sty   
\magnification 1200
\vsize = 9.5 true in
\hsize=6.2 true in
\NoRunningHeads        
\parskip=\medskipamount
        \lineskip=2pt\baselineskip=18pt\lineskiplimit=0pt
       
        \TagsOnRight
        \NoBlackBoxes

        \topmatter
        \title
        Logarithmic Bounds on Sobolev Norms for Time Dependent Linear Schr\"odinger Equations      
        \endtitle
        \author
        W.-M.~Wang
        \endauthor
\address
Departement de Mathematique, Universite Paris Sud, 91405 Orsay Cedex, FRANCE
\endaddress
        \email
{wei-min.wang\@math.u-psud.fr}
\endemail
\thanks  
\endthanks

        \bigskip\bigskip
        \bigskip
        \toc
        \bigskip
        \bigskip 
        \widestnumber\head {Table of Contents}
        \head 1. Introduction and statement of the theorem
        \endhead
        \head 2. Periodic approximations and Floquet solutions 
        \endhead
        \head 3. Some a priori estimates 
        \endhead
        \head 4. Bounds on Sobolev norms
        \endhead
        \endtoc
        \endtopmatter
        \vfill\eject
        \bigskip
\document
\head{\bf 1. Introduction and statement of the theorem}\endhead
We consider the time dependent linear Schr\"odinger equation: 
$$i\frac\partial{\partial t}u=-\Delta u+V(x, t)u,
\tag 1.1$$
on $\Bbb T{\overset\text{def }\to=}[-\pi, \pi)$ 
with periodic boundary conditions. The potential $V$ is identified with a function on $\Bbb R\times\Bbb R$, periodic in $\Bbb R$ with period $2\pi$. (To emphasize the time dependence, we write $V(x, t)$
for $V$.) We further assume that $V$ is real analytic in $(x, t)$ in a strip 
$D{\overset\text{def }\to=}(\Bbb R+i\rho)^2$ ($|\rho|<\rho_0$, $\rho_0>0$), real and bounded 
in $\Bbb R^2$: $\Vert V\Vert_{\infty, \Bbb R^2}<C<\infty$. We prove the following result:
\proclaim{Theorem} There exists $\varsigma>3$, such that for all $s>0$, there exists $C_s$, such that
$$\Vert u(t)\Vert_{H^s}\leq C_s[\log (|t|+2)]^{\varsigma s}\Vert u(0)\Vert_{H^s},\tag 1.2$$
where $u(t)$ is the solution to (1.1) with the initial condition $u_0\in H^s$.
\endproclaim

\noindent{\it Remark.} As we will see later, the theorem in fact extends to the class of potentials
where the analyticity is replaced by some uniform estimates on the derivatives (e.g., Gevrey) and
boundedness on $\Bbb R^2$ is weakened to logarithmic growth in $t$. The exponent $\varsigma$
will however depend on the Gevrey exponent. In this paper, we only state the theorem in the 
analytic case. Previously, it was proven in \cite{B2, 3} that for $V\in C^\infty$,
$$\Vert u(t)\Vert_{H^s}\leq C_{s, \epsilon}(|t|+1)^{\epsilon}\Vert u(0)\Vert_{H^s},$$
for all $\epsilon>0$.

In \cite{N}, it was proven that for smooth time dependent potentials with certain random dependence on time,
$\Vert u(t)\Vert_{H^s}$ is almost surely unbounded in time, which shows that the $\log$ in (1.2) is almost
surely necessary. On the other hand, in \cite{W}, it was proven that for an explicit time periodic potential 
$\Vert u(t)\Vert_{H^s}$ remains bounded for all $t$. Clearly this belongs to the exceptional set from the
point of view of random dependence in time.  The present Theorem together with \cite{N, W}
give a rather complete picture of time dependent linear Schr\"odinger equations on the circle. 

The proof consists of making periodic in time approximations by replacing $V(x, t)$ with $V_1(x, t)$
which is periodic in $t$ with period $2\pi T$ and $V_1(x, t)=V(x, t)$ for $|t|\leq T$. The dynamics of 
equation (1.1) is hence equivalent to the dynamics of 
$$i\frac\partial{\partial t}u=-\Delta u+V_1(x, t)u, \tag 1.3$$
for $|t|\leq T$. This part of the strategy is similar to \cite{B2, 3}.

From Floquet theory, the dynamics of (1.3) can be reduced to the spectral theory of the corresponding
Floquet operator
$$H=\text { diag }(\frac{n}{T}+j^2)+\hat V_1*\tag 1.4$$
on $\ell^2 (\Bbb Z^2)$, where $n$ is the dual variable of $t$, $j$ the dual of $x$ and
$$\hat V_1(j, n)=\int_{-\pi T}^{\pi T}\int_{-\pi }^{\pi}V_1(x, t)e^{-ijx}e^{-i\frac{n}{T}t}dxdt\tag 1.5$$
is the Fourier transform of $V_1$.

More specifically, let $\Bbb T_T=[-\pi T, \pi T)$ with periodic boundary conditions. For any initial
datum $u_0\in L^2(\Bbb T)$, we identify $u_0$ with $\check u_0\in L^2(\Bbb T)\times L^2(\Bbb T_T) $
as follows:
$$
\cases
\tilde u_0(j, 0)=\hat u_0(j),\\
\tilde u_0(j, n)=0,\, n\neq 0\endcases\tag 1.6
$$
where $\tilde u_0$ is the Fourier transform of $\check u_0$ and $\hat u_0$ that of $u_0$.  All solutions to (1.3) can be written as linear superpositions of Floquet solutions, which up to a phase are inverse Fourier transforms of eigenfunctions to (1.4). Localization of eigenfunctions of (1.4) therefore leads
to control over Sobolev norms of solutions to (1.3).

When $T$ is a {\it fixed} integer period, which is the resonant case as $\sigma(\Delta)=\{j^2|j\in\Bbb Z\}\subset\Bbb N$, it was proven in \cite{W} that under appropriate conditions (1.4) has pure point
spectrum with exponentially localized eigenfunctions. This in turn leads to
$$\Vert u(t)\Vert_{H^s}\leq C_s\Vert u(0)\Vert_{H^s}$$ for all $s>0$.

The main complication here is that $T$ is a variable, in fact $T\to\infty$. Using the identification (1.6) 
and the fact that for $|t|\leq T$ and large frequencies $|j|>J(T)$, the $H^s$ norm is essentially 
preserved, see (3.5) of Lemma 3.1, we construct approximate eigenfunctions to (1.4) by 
restricting $H$ to $H_\Lambda$ with $\Lambda=\{(j, n)||j|\leq J(T),\, |n|\leq N(T)\}$ for approprate
$J$, $N$ depending on $T$. This differs from \cite{B2, 3} and enables us to obtain logarithmic bounds.

We prove that these approximate eigenfunctions are localized in the proposition in sect. 2. The
separation properties of the set $\{j^2|j\in\Bbb Z\}$ again plays an important role here as in \cite{W}.
In sect. 4, using the proposition and some a priori estimate for linear Schr\"odinger equations, we
prove the theorem.
\smallskip
\head{\bf 2. Periodic approximations and Floquet solutions}\endhead

Let $u(t)$ be the solution to 
$$i\frac\partial{\partial t}u=-\Delta u+V(x, t)u,\tag 2.1$$
with the initial condition $u_0\in H^s$ ($0<s<\infty$). We want to 
bound $\Vert u(t)\Vert_{H_s}$ as $t\to \infty$. (When $s=0$, the $L^2$ norm is conserved.)
We therefore look at (2.1) for $0\leq |t|\leq \pi T$, $T\gg 1$, with the initial condition
$u_0$.  Let $\tilde\phi\in C_0^\infty[-\pi, \pi]$ be a (fixed) Gevrey function of order $\alpha$:
$$\max_{\tau\in[-\pi,\pi]}\big|\frac{\partial^m\tilde\phi(\tau)}{\partial\tau^m}\big|\leq C^{m+1}(m!)^\alpha,\qquad 1<\alpha<\infty\tag 2.2$$
satisfying 
$$
\cases 0\leq\tilde\phi\leq 1,\\
\tilde\phi(\tau)=1,\, |\tau|\leq 1\\
\tilde\phi(\tau)=0,\, |\tau|\geq \pi,\endcases\tag 2.3
$$
(cf. \cite{H\"o}).

Let $$\phi(t)=\tilde\phi(\frac{t}{T}).\tag 2.4$$
Define
$$V_1(x, t)=\sum_{j\in\Bbb Z}V(x, t+2\pi jT)\phi(t+2\pi j T).\tag 2.5$$
Then $V_1(x, t)$ is $2\pi$ periodic in $x$, $2\pi T$ periodic in $t$, analytic in $x$, Gevrey in $t$ 
of order $\alpha$, ($1<\alpha<\infty$). 
$$V_1(x, t)=V(x, t)\tag 2.6$$
for  $0\leq |t|\leq T$, and $$\Vert V_1\Vert_\infty\leq 2 \Vert V\Vert_\infty.\tag 2.7$$
So for $0\leq |t|\leq T$, we can study instead the equation
$$i\frac\partial{\partial t}u=-\Delta u+V_1(x, t)u.
\tag 2.8$$
$V_1$ has the Fourier decomposition:
$$V_1(x, t)=\sum_{j\in\Bbb Z,\, n\in\Bbb Z}\hat V_1(j, n) e^{i(jx+\frac {n}{T} t)},\tag 2.9$$
where $$\align|\hat V_1(j, n)|&\leq C e^{-c|j|},\qquad \quad |j|\geq (\log T)^\delta,\\
&\leq C e^{-c|\frac{n}{T}|^{1/\alpha}},\quad |n|\geq T(\log T)^\delta,\tag 2.10\\
&\qquad (T\gg 1, \, 0<C, c,\delta<\infty).\endalign$$

Since we seek solutions to (2.8) for finite time: $|t|\leq T$, it is convenient to replace $V_1$ by $V_2$ 
defined as 
$$V_2(x, t)=\sum_{\Sb |j|\leq (\log T)^{\sigma}\\ |n|\leq T(\log T)^{\sigma}\endSb}
\hat V_1(j, n) e^{i(jx+\frac {n}{T} t)},\tag 2.11$$
where $\varsigma>\sigma>\alpha+\delta>1$. Using (2.10), 
$$\Vert V_1-V_2\Vert_\infty\leq e^{-(\log T)^{\sigma'/\alpha}}\ll \frac{1}{T^p}\tag 2.12$$
provided $p<(\log T)^{\sigma'/\alpha-1}$ ($T\gg 1$), where $1<\alpha+\delta<\sigma'<\sigma$. 
For $|t|\leq T$, (2.12) will permit us to use Floquet solutions to 
$$i\frac\partial{\partial t}u=-\Delta u+V_2(x, t)u,
\tag 2.13$$
in the approximation process in view of the following basic fact.

\proclaim{Lemma 2.1} Let $\tilde u$ be an approximative solution of (2.8):
$$(i\frac\partial{\partial t}+\Delta -V(x, t))\tilde u=\eta,$$
with $\tilde u(t=0)=u_0$, where $\Vert\eta(t)\Vert_{L^2}\leq\epsilon$ for all $|t|\leq T$. Then the solution
$u$ to (2.8) with $\tilde u(t=0)=u_0$ satisfies 
$$\Vert u(t)-\tilde u(t)\Vert_{L^2}<\epsilon|t|\leq\epsilon T$$ for $|t|\leq T$.
\endproclaim

\demo{Proof}
Let $S(t)$ denote the flow of (2.8). This follows from the integral equation
$$(\tilde u-u)(t)=i\int_0^t S(t)S(\tau)^{-1}\eta(\tau)d\tau$$ and $\Vert S(t)\Vert_{L^2\to L^2}=1$. $\square$
\enddemo
\smallskip
\noindent{\it Floquet solutions to (2.13)}

Since (2.13) is time periodic with period $2\pi T$, any $L^2$ solution can be written as a linear 
superposition of Floquet solutions of the form $e^{iEt}\psi(x, t)$, where $\psi(x, t)$ is $2\pi$ periodic 
in $x$ and $2\pi T$ periodic in $t$:
$$\psi (x, t)=\sum_{(j, n)\in\Bbb Z^2}
\hat \psi(j, n) e^{i(jx+\frac {n}{T} t)},\tag 2.14$$
$E$ is called the Floquet eigenvalue; $E$, $\hat\psi$ satisfy the eigenvalue equation:
$$\aligned H\hat\psi&=[\text {diag }(\frac{n}{T}+j^2)+\hat V_2*]\hat\psi\\
&=E\hat\psi\endaligned\tag 2.15$$
on $\ell^2(\Bbb Z^2)$, where $*$ denotes convolution:
$$(\hat V_2*\hat\psi)(j, n)=\sum_{(j', n')\in\Bbb Z^2}\hat V_2(j-j', n-n')\hat\psi(j', n'),\tag 2.16$$
$$\aligned \hat V_2(j, n)&=\hat V_1(j, n)\qquad\text{if }|j|\leq (\log T)^{\sigma}\text{ and }|n|\leq T(\log T)^{\sigma},\, \sigma>\alpha+\delta>1,\\
&=0\qquad\qquad\quad\text{otherwise}, \endaligned\tag 2.17$$
and $\hat V_1$ satisfies (2.10).

We identify the initial condition $\hat u_0\in\ell^2(\Bbb Z)$ with $\tilde u_0\in\ell^2(\Bbb Z^2)$, where 
$$\cases \tilde u_0(j, 0)=\hat u_0(j)\\ \tilde u_0(j, n)=0,\quad n\neq 0.\endcases\tag 2.18$$
Since we are only concerned about finite time: $|t|\leq T$, in view of (2.12, 2.17, 2.18), Lemma 2.1, the a priori estimate (3.5)
and some related estimates, which we will prove in sect. 3 (see Lemma 3.1), it is sufficient to solve the eigenvalue problem in (2.15) in a finite region
$$\Lambda=\{(j, n)\in\Bbb Z^2|\,|j|\leq J(T),\, |n|\leq AT(\log T)^\sigma\},\tag 2.19$$ 
where $J(T)>T^s$ depending on $T$ and the Sobolev index $s$, $A>1$ as in the following proposition, $\sigma>\alpha+\delta>1$ as in (2.1). 

For any subset $\Cal S\subset\Bbb Z^2$, define $H_{\Cal S}$ to be the restriction of $H$ to $\Cal S$:
$$H_{\Cal S}(n,j;n'j')=\cases H(n,j;n'j'),\qquad (n,j)\text{ and }(n',j')\in\Cal S\\0\quad\,\qquad\qquad\qquad\text{otherwise}.
\endcases\tag 2.20$$
We have the following estimates on eigenfunctions of $H_\Lambda$.
\proclaim{Proposition} Assume $$H_\Lambda\xi=E\xi,\qquad\Vert \xi\Vert_{\ell^2(\Lambda)}=1.\tag 2.21$$
Define $$\Omega_0=\{(j, n)\in\Lambda|\, |j|\leq 4A(\log T)^\sigma\},\quad (\sigma>\alpha+\delta>1)\tag 2.22$$
and for any $(j_0, n_0)\in\Lambda$, define
$$\Omega'(j_0, n_0)=\{(j, n)\in\Lambda|\,|\, |j|-|j_0|\,|\leq (\log T)^\sigma, \, |n-n_0|\leq T(\log T)^\sigma\},\,
(\sigma>\alpha+\delta>1).\tag 2.23$$
Then for all $\xi$ eigenfunctions of $H_\Lambda$ as in (2.21), $\xi$ satisfies either 
$$\align\Vert \xi\Vert_{\ell^2(\Lambda\backslash\Omega_0)}&\leq e^{-(\log T)^{(\frac{\sigma'-\delta}{\alpha})}}\tag 2.24\\
\text{or }\Vert \xi\Vert_{\ell^2(\Lambda\backslash\Omega')}&\leq e^{-(\log T)^{(\frac{\sigma'-\delta}{\alpha})}},\,
(1<\alpha+\delta<\sigma'<\sigma).\tag 2.25\endalign$$
for some $\Omega'=\Omega'(j_0, n_0)$, $(j_0, n_0)\in\Lambda$.
\endproclaim
\demo{Proof} For any given $E$, we define the resonant set $\Omega$ such that if $(j, n)\in\Omega$, then 
$$|\frac{n}{T}+j^2-E|\leq (\log T)^\sigma,\quad (\sigma>\alpha+\delta>1).\tag 2.26$$
So $$\Vert (H_{\Lambda\backslash\Omega}-E)^{-1}\Vert \leq\frac{1}{(\log T)^\sigma-\Vert V_2\Vert_\infty}\leq \frac{2}{(\log T)^\sigma},\tag 2.27$$
if  $$\Vert V_2\Vert_\infty<\frac{1}{2}(\log T)^\sigma.$$
From (2.26), we have 
$$E-\frac{n}{T}-(\log T)^\sigma\leq j^2\leq E-\frac{n}{T}+(\log T)^\sigma,\quad (\sigma>\alpha+\delta>1)\tag 2.28$$
for $(j, n)\in\Omega$. We distinguish the following two cases:
\item{(i)} $E\leq 5A^2 (\log T)^{2\sigma}$

The less or equal part of (2.28) gives $ |j|\leq 3A(\log T)^\sigma$. So $\Omega\subset\{(j,n)\in\Lambda||j|\leq 
3A(\log T)^\sigma\}$. Define $B=\Lambda\backslash\Omega$, $B_0=\Lambda\backslash\Omega_0$, $B_0\subset B$. Let $P_B$, $P_{B_0}$ be projections onto the sets $B$, $B_0$. 

Assume $\xi$ is an eigenfunction with eigenvalue $E\leq 5A^2 (\log T)^{2\sigma}$. Then
$$P_B\xi=-(H_B-E)^{-1}P_B\Gamma\xi\tag 2.29$$
where $$\Gamma=H_\Lambda-H_B\oplus H_\Omega.\tag 2.30$$
So $$P_B\xi=-(H_B-E)^{-1}P_B\Gamma P_\Omega\xi.\tag 2.31$$
Let $$\Gamma_0=H_B-H_{B_0}\oplus H_{B\backslash B_0}.\tag 2.32$$ 
Then 
$$\aligned P_{B_0}\xi&=P_{B_0}P_B\xi\\
&=-P_{B_0}(H_{B_0}-E)^{-1}P_{B_0}\Gamma P_\Omega\xi\\
&\quad +P_{B_0}(H_{B_0}-E)^{-1}\Gamma_0(H_{B}-E)^{-1}P_{B}\Gamma P_\Omega\xi,\endaligned\tag 2.33$$
where we used $B_0\subset B$. Using (2.27) on $(H_{B_0}-E)^{-1}$ and $(H_{B}-E)^{-1}$  and (2.10, 2.11, 2.15)
on $\Gamma$, $\Gamma_0$, we obtain 
$$\Vert P_{B_0}\xi\Vert_{\ell^2}\leq \frac{4e^{-c(\log T)^{\sigma-\delta}}}
{(\log T)^{2\sigma}}<
e^{-(\log T)^{\frac{\sigma'-\delta}{\alpha}}}\quad (1<\alpha+\delta<\sigma'<\sigma,\, T\gg 1),\tag 2.34$$
which is (2.24).

\item{(ii)} $E> 5A^2 (\log T)^{2\sigma}$

The greater or equal part of (2.28) gives 
$$|j|\geq 2A(\log T)^\sigma.\tag 2.35$$
So if there exist $(j, n)$, $(j', n')\in\Omega\subset\Lambda$, $|j|\neq |j'|$, then 
$$|\frac{n-n'}{T}+j^2-{j'}^2|\leq 2 (\log T)^\sigma$$
from (2.26). Using (2.35), this implies
$$\aligned |\frac{n-n'}{T}|&\geq (|j|+|j'|)(|j|-|j'|)-2 (\log T)^\sigma\\
&\geq (4A-2)(\log T)^\sigma\\
&>2A(\log T)^\sigma\endaligned$$
if $A>1$, which is a contradiction from the definition of $\Lambda$. So $|j|=|j'|$ and 
$$|\frac{n-n'}{T}|\leq 2(\log T)^\sigma<2A(\log T)^\sigma$$
for $A>1$, if both $(j, n)$, $(j', n')\in\Omega$. (2.25) follows by using the same argument as in (2.29-2.34)
with $\Omega'$ replacing $\Omega_0$. $\square$
\enddemo
\smallskip
\head{\bf 3. Some a priori estimates}\endhead

In this section, we collect some basic estimates on the flow of linear Schr\"odinger equations with smooth
potentials, cf. \cite{B1}. Since we will need estimates on $H^s$ norms for $s>0$ dependent on $T$, for 
completeness we also include their proofs, making explicit the dependence on $s$. 

Let $S(t)$ be the flow of the linear Schr\"odinger equation in (2.1). Then $S(t)$ is unitary: $\Vert S(t)\Vert_{\ell^2\to\ell^2}=1$. Let $\Pi_J$ denote the Fourier multiplier defined as 
$$\aligned\hat \Pi_J&=1,\qquad\qquad\qquad\qquad|j|\leq J/2,\\
&=2(1-|j|/J),\qquad\quad\, J/2\leq|j|\leq J,\\
&=0,\qquad\qquad\qquad\qquad |j|>J.\endaligned\tag 3.1$$
Since $V$ is real analytic in $(x, t)$ and bounded in $D$, we have 
$$\Vert \frac{\partial^m V}{\partial x^m}\Vert_{\infty,\,\Bbb T}\leq C^{m+1} m!,\qquad m=0,\, 1,...\tag 3.2$$
We have the following estimates on the $H^s$ norms:
\proclaim{Lemma 3.1}
$$\align &\Vert S(t)\Vert_{H^s\to H^s}\leq C^s s! (|t|^s+1),\tag 3.3\\
&\Vert [V,\, \Pi_J]\Vert_{H^s\to H^s}\leq \frac{Cs!}{J}\quad (J\gg 1),\tag 3.4\\
&\Vert (I-\Pi_J)S(t)\Vert_{H^s\to H^s}\leq 1+\frac{(C^s s!)^2 }{J} |t|^{s+1}\quad (J>|t|^s),\tag 3.5\\
&\Vert [S(t), \Pi_J]\Vert_{H^s\to H^s}\leq \frac{(C^s s!)^4}{J} (|t|^{3s+1}+1)\quad (J>|t|^s).\tag 3.6\endalign$$\endproclaim
\noindent {\it Remark.} The same estimates hold for the flow of (2.8) as only the $x$-derivatives are involved.
\demo{Proof} Using (2.1), 
$$\aligned \frac{\partial}{\partial t}\Vert& u(t)\Vert^2_{H^s}=2\text{ Re }(\frac{\partial^s}{\partial x^s}u(t), 
\frac{\partial^s}{\partial x^s}\frac{\partial}{\partial t}u(t))\\
=&2\text{ Im }(\frac{\partial^s}{\partial x^s}u(t), 
\frac{\partial^s}{\partial x^s}(\Delta u+Vu))\\
=&2\text{ Im }(\frac{\partial^s}{\partial x^s}u(t), 
\sum_{\Sb \gamma+\beta=s\\\gamma\geq 1\endSb}\frac{\partial^\gamma V}{\partial x^\gamma}\frac{\partial^\beta u}{\partial x^\beta}).\endaligned\tag 3.7$$
It follows that 
$$\frac{\partial}{\partial t}\Vert u(t)\Vert_{H^s}\leq \sum_{\Sb \gamma+\beta=s\\\gamma\geq 1\endSb}
\Vert u(t)\Vert_{H^\beta}C^{\gamma+1}\gamma!,\tag 3.8$$
where we used (3.2).

Using interpolation:
$$\Vert u(t)\Vert_{H^{s-1}}\leq \Vert u(t)\Vert_{H^{s}}^{\frac{s-1}{s}}\Vert u(t)\Vert_{L^2}^{\frac{1}{s}} \quad (s\geq 1),$$
and more generally,
$$\Vert u(t)\Vert_{H^{s-\gamma}}\leq \Vert u(t)\Vert_{H^{s}}^{\frac{s-\gamma}{s}}\Vert u(t)\Vert_{L^2}^{\frac{\gamma}{s}} \quad (s\geq \gamma).\tag 3.9$$
Using (3.9) in (3.8), we have
$$\aligned \frac{\partial}{\partial t}\Vert u(t)&\Vert_{H^{s}}\leq C^2\Vert u(t)\Vert_{H^{s}}^{1-\frac{1}{s}}\Vert u(t)\Vert_{L^2}^{\frac{1}{s}}+ C^3\Vert u(t)\Vert_{H^{s}}^{1-\frac{2}{s}}\Vert u(t)\Vert_{L^2}^{\frac{2}{s}}+\cdots\\
&+ C^{\gamma+1}\gamma!\Vert u(t)\Vert_{H^{s}}^{1-\frac{\gamma}{s}}\Vert u(t)\Vert_{L^2}^{\frac{\gamma}{s}}+\cdots+C^{s+1}s!\Vert u(t)\Vert_{L^2}\endaligned\tag 3.10$$
Since $$\frac{s}{\gamma}\Vert u(t)\Vert_{H^{s}}^{1-\frac{\gamma}{s}}\frac{\partial}{\partial t}\Vert u(t)\Vert_{H^s}^{\frac{\gamma}{s}}=\frac{\partial}{\partial t}\Vert u(t)\Vert_{H^s}\quad (1\leq\gamma\leq s),\tag 3.11$$
we obtain from (3.10),
$$\Vert u(t)\Vert_{H^{s}}\leq C^s s! (|t|^s+1)\Vert u_0\Vert_{H^{s}}.\tag 3.12$$
Hence
$$\Vert S(t)\Vert_{H^s\to H^s}\leq C^s s! (|t|^s+1).\tag 3.13$$

To prove (3.4), it is more convenient to work with the Fourier variables $j$ dual to $x$. Let ${\hat {\hat V}}$
be the partial Fourier transform with respect to $x$. we have 
$$[V, \Pi_J]\hat(j, j')={\hat {\hat V}}(j-j')(\hat\Pi_J(j')-\hat\Pi_J(j)),\tag 3.14$$
where $\hat\Pi_J$ is defined in (3.1). Since $V$ is analytic, periodic in $x$ and $|V(x, t)|<C$ for all $t$,
$$\aligned |{\hat {\hat V}}(j-j')|&\leq C e^{-c|j-j'|},\\
\text{and from }(3.1)\qquad|\hat\Pi_J(j')-\hat\Pi_J(j)|&\leq 1,\qquad\qquad\quad  |j-j'|\geq J/2,\\
&\leq\frac{2}{J}|j-j'|,\quad\,\, |j-j'|<J/2.\endaligned\tag 3.15$$ 
Using (3.15), we have 
$$\aligned |[V, \Pi_J]\hat(j, j')|&\leq C e^{-c|j-j'|}, \quad\quad\quad\quad\quad|j-j'|\geq J/2\\
&\leq \frac{2C}{J}|j-j'|e^{-c|j-j'|},\,\quad|j-j'|< J/2\endaligned \tag 3.16$$
From Schur's lemma, we then obtain (3.4).

To prove (3.5), we proceed similarly to the proof of (3.3). We have 
$$\aligned \frac{\partial}{\partial t}\Vert (I-\Pi_J)u(t)\Vert_{H^s}^2&=2\text{ Im }((I-\Pi_J)\frac{\partial^s}{\partial x^s}u(t), (I-\Pi_J)V\frac{\partial^s}{\partial x^s}u(t))\\
&\quad +2\text{ Im }((I-\Pi_J)\frac{\partial^s}{\partial x^s}u(t), (I-\Pi_J)\sum_{\Sb \gamma+\beta=s\\\gamma\geq 1\endSb}\frac{\partial^\gamma V}{\partial x^\gamma}\frac{\partial^\beta u}{\partial x^\beta}).\endaligned\tag 3.17$$
So
$$\aligned  &\frac{\partial}{\partial t}\Vert (I-\Pi_J)u(t)\Vert_{H^s}\\
\leq &[V,\Pi_J]\Vert u(t)\Vert_{H^s}+[\frac{\partial V}{\partial x},\Pi_J]\Vert u(t)\Vert_{H^{s-1}}+C\Vert(I-\Pi_J)u(t)\Vert_{H^{s-1}}\\
&+\cdots+[\frac{\partial^\gamma V}{\partial x^\gamma},\Pi_J]\Vert u(t)\Vert_{H^{s-\gamma}}+C^{\gamma+1}\gamma!\Vert(I-\Pi_J)u(t)\Vert_{H^{s-\gamma}}\\
&+\cdots+[\frac{\partial^s V}{\partial x^s},\Pi_J]\Vert u(t)\Vert_{L^2}+C^{s+1}s!\Vert(I-\Pi_J)u(t)\Vert_{L^2}\endaligned\tag 3.18$$
Using $$\align &\Vert [\frac{\partial^\gamma V}{\partial x^\gamma},\Pi_J]\Vert_{H^{s-\gamma}\to H^{s-\gamma}}\leq\frac{Cs!}{J},\tag 3.19\\
&\Vert (I-\Pi_J)u(t)\Vert_{H^{s-\gamma}}\leq \frac{1}{J^{\gamma}}\Vert u(t)\Vert_{H^s}\tag 3.20\endalign$$
and (3.9) in (3.18) and integrating over $t$, we obtain (3.5).

To prove (3.6), assume $u$ is a solution to (2.1)
$$i\frac\partial{\partial t}u+\Delta u-V(x, t)u=0,$$
then $$(i\frac\partial{\partial t}+\Delta)\Pi_J u-V(\Pi_J u)=-[V, \Pi_J] u.\tag 3.21$$
From Lemma 2.1
$$\aligned &[S(t), \Pi_J] u_0\\
=&S(t)\Pi_J u_0-\Pi_J S(t)u_0\\
=&-i\int_0^t S(t)S(\tau)^{-1}[V, \Pi_J ]u(\tau)d\tau\endaligned\tag 3.22$$
Using (3.3, 3.4) in (3.22), we obtain (3.6). $\square$
\enddemo
\smallskip
\head{\bf 4. Bounds on Sobolev norms}\endhead

Let $u_0\in H^s$ be an initial datum, normalized so that $\Vert u_0\Vert_{H^s}=1$. We assume $0<s\leq\log T$,
cf. (4.49). Let $$J=T^{10s}.\tag 4.1$$ Then from (3.5)
$$\aligned \Vert S(0, t)u_0\Vert_{H^s}&=\Vert \Pi_{J/4}S(0, t)u_0\Vert_{H^s}+\Vert (1-\Pi_{J/4})S(0, t)u_0\Vert_{H^s}\\
&\leq \Vert \Pi_{J/4}S(0, t)u_0\Vert_{H^s}+2\Vert u_0\Vert_{H^s}.\endaligned\tag 4.2$$

In view of the proposition and (4.2), let $$J_0=4A(\log T)^\sigma\quad (A>1, \,\sigma>\alpha+\delta>1).\tag 4.3$$
We make the following decomposition:
$$\align \Vert \Pi_{J/4}S(0, t)u_0\Vert_{H^s}&\leq  \Vert \Pi_{J/4}S(0, t)\Pi_{2J_0}u_0\Vert_{H^s}\tag 4.4\\
&+\Vert \Pi_{J/4}S(0, t)(\Pi_{J/2}-\Pi_{2J_0})u_0\Vert_{H^s}\tag 4.5\\
&+\Vert \Pi_{J/4}S(0, t)(I-\Pi_{J/2})u_0\Vert_{H^s}.\tag 4.6\endalign$$
(4.6) can be estimated using (3.5):
$$\aligned \Vert \Pi_{J/4}S(0, t)(I-\Pi_{J/2})u_0\Vert_{H^s}&\leq \Vert [S(0, t),\Pi_{J/2}]u_0\Vert_{H^s}\\
&\leq 2\frac{(C^s s!)^4}{J}(1+T^{3s+1})\Vert u_0\Vert_{H^s}\\
&<1.\endaligned\tag 4.7$$

To estimate (4.5), we use the following
\proclaim{Lemma 4.1} Let $\phi$ be such that $$\text {supp }\hat\phi\subseteq[-J/2, -2J_0]\cup[2J_0, J/2],
\tag 4.8$$ then 
$$\Vert \Pi_{J/4}S(t)\phi\Vert_{H^s}\leq C^s\Vert\phi\Vert_{H^s}.\tag 4.9$$
\endproclaim
\demo{Proof} We identify $\hat\phi$ with $\tilde\phi$ defined as 
$$\cases\tilde\phi(j, 0)=\hat\phi(j),\\ \tilde\phi(j, n)=0,\quad n\neq 0.\endcases\tag 4.10$$
$\text{supp }\tilde\phi\subset\Lambda$, where $\Lambda$ is defined in (2.19). $\tilde\phi\in\ell^2(\Lambda)$.
So we can expand $\tilde\phi$ using the eigenfunctions $\xi$ of $H_\Lambda$:
$$\tilde\phi=\sum(\tilde\phi,\xi)\xi.\tag 4.11$$ 

Let $\chi_S$ be the characteristic function of the set $S$: 
$$\chi_S|_S=1,\quad \chi_S|_{\Lambda\backslash S}=0.\tag 4.12$$
For an eigenfunction $\xi$ satisfying (2.25), let 
$$\xi'=\chi_{\Omega'}\xi,\quad Q=\{\xi'|\xi\text{ satisfies } (2.25)\}.\tag 4.13$$
Using (4.8, 4.10, 2.24, 2.25) in (4.11), we have 
$$\Vert \tilde\phi-\sum_{\xi'\in Q}(\tilde\phi,\xi')\xi'\Vert_{\ell^2(\Lambda)}\leq \Cal O(e^{-(\log T)^{(\frac{\sigma'-\delta}{\alpha})}}\sqrt {|\Lambda|}\Vert\hat\phi\Vert_{\ell^2}).\tag 4.14$$
Since $|\Lambda|\leq T^{10s+2}$ from (2.19, 4.1) and $0<s\leq\log T$, we have 
$$\Vert \tilde\phi-\sum_{\xi'\in Q}(\tilde\phi,\xi')\xi'\Vert_{\ell^2(\Lambda)}\leq e^{-\frac{2}{3}(\log T)^{(\frac{\sigma'-\delta}{\alpha})}}, \tag 4.15$$ 
assuming $$\sigma'>2\alpha+\delta>2.\tag 4.16$$

From (4.13, 2.25), $\xi'$ is an approximate eigenfunction of $H_\Lambda$:
$$\aligned\Vert (H_\Lambda-E)\xi'\Vert_{\ell^2(\Lambda)}&\leq\Vert (H_\Lambda-E)\Vert_{\ell^2\to\ell^2}
\Vert\xi\Vert_{\ell^2(\Lambda\backslash\Omega')}\\
&=e^{-\frac{2}{3}(\log T)^{(\frac{\sigma'-\delta}{\alpha})}}, \quad (\sigma'>2\alpha+\delta>2).\endaligned\tag 4.17$$
Hence $\xi'$ is an approximate eigenfunction of $\tilde H$:
$$\tilde H {\overset\text{def }\to=}\text{ diag }(j^2+\frac{n}{T})+\hat V_1*=H+(\hat V_1-\hat V_2)*,\tag 4.18$$
where $H$ as defined in (2.15), $\hat V_1$ is defined in (2.5, 2.9) satisfying (2.10).  This is because
$$(\tilde H-E)\xi'=(H_\Lambda-E)\xi'+\Gamma\xi'+(\hat V_2-\hat V_1)*\xi',\tag 4.19$$
where $H_\Lambda$ as defined in (2.20, 2.19), $$\Gamma=H-H_\Lambda\oplus H_{\Lambda^c}\tag 4.20$$
and 
$$\Vert (\tilde H-E)\xi'\Vert_{\ell^2}\leq 2 e^{-\frac{2}{3}(\log T)^{(\frac{\sigma'-\delta}{\alpha})}}\quad (\sigma'>2\alpha+\delta>2),\tag 4.21$$
using (4.17, 2.12).

Define $$\check\xi(x, t)=e^{iEt}\sum_{(j,n)\in\Omega'}\xi'(j,n) e^{i(jx+\frac{n}{T}t)}.\tag 4.22$$
From (4.18), $\check\xi$ is an approximate Floquet solution of (2.8) satisfying 
$$(i\frac{\partial}{\partial t}+\Delta -V_1)\check\xi=e^{-\frac{2}{3}(\log T)^{(\frac{\sigma'-\delta}{\alpha})}}\quad 
(\sigma'>2\alpha+\delta>2).
\tag 4.23$$
Let $S(t)$ be the flow for equation (2.8), using Lemma 2.1
$$\Vert \check\xi(t)-S(t)\check\xi(0)\Vert_2\leq Te^{-\frac{2}{3}(\log T)^{(\frac{\sigma'-\delta}{\alpha})}}.\tag 4.24$$
Because of localization properties of $\xi'$ in (2.25), we will have good control over $\Vert \check\xi(t)\Vert_{H^s}$
(see (4.28-4.37)).

In view of (4.24), we express $\phi$ as an approximate linear combination of $\check\xi(0)$ as follows. Taking
the inverse Fourier transform of the expression under the norm sign in (4.14), we have equivalently
$$\Vert \tilde\phi(x,\theta)-\sum_{\xi'\in Q}(\tilde\phi,\xi')\sum_{(j,n)\in\Omega'}\xi'(j,n)e^{i(jx+\frac{n}{T}\theta)}\Vert_
{L^2(\Bbb T\times\Bbb T_T)}
\leq e^{-\frac{2}{3}(\log T)^{(\frac{\sigma'-\delta}{\alpha})}},\tag 4.25$$
where $\Bbb T$ denotes $[-\pi, \pi)$ with periodic boundary conditions and $\Bbb T_T$ denotes $[-\pi T, \pi T)$ with periodic boundary conditions.

So $\phi(x){\overset\text{def }\to=}\phi(x, 0)$ as a function on $L^2(\Bbb T)$ satisfies 
$$\aligned&\Vert \phi(x)-\sum_{\xi'\in Q}(\tilde\phi,\xi')\sum_{(j,n)\in\Omega'}\xi'(j,n)e^{ijx}\Vert_{L^2(\Bbb T)}\\
\leq &T^{1/2}(\log T)^{\sigma/2}e^{-\frac{2}{3}(\log T)^{(\frac{\sigma'-\delta}{\alpha})}}\\
\leq &e^{-\frac{1}{2}(\log T)^{(\frac{\sigma'-\delta}{\alpha})}},\,(\sigma>\sigma'>2\alpha+\delta>2).\endaligned\tag 4.26$$
Therefore for $|t|\leq T$, (4.24, 4.26) give 
$$\aligned&\Vert S(t)\phi-\sum_{\xi'\in Q}(\tilde\phi,\xi')\check\xi(t)\Vert_{L^2(\Bbb T)}\\
\leq &e^{-\frac{1}{2}(\log T)^{(\frac{\sigma'-\delta}{\alpha})}}+\sum_{\xi'\in Q}|(\tilde\phi,\xi')\Vert \check\xi(t)-S(t)\check\xi(0)\Vert_{L^2(\Bbb T)}\\
\leq &e^{-\frac{1}{2}(\log T)^{(\frac{\sigma'-\delta}{\alpha})}}+|\Lambda|Te^{-\frac{2}{3}(\log T)^{(\frac{\sigma'-\delta}{\alpha})}}\\
\leq &2e^{-\frac{1}{2}(\log T)^{(\frac{\sigma'-\delta}{\alpha})}}\quad (\sigma'>2\alpha+\delta>2),\endaligned
\tag 4.27$$
where we used (4.1, 2.19).

So we only need to estimate $\Vert \sum_{\xi'\in Q}(\tilde\phi,\xi')\check\xi(t)\Vert_{H^s}$ ($s>0$). Let $\delta_j$ be the Dirac delta function at $j$. We have 
$$\aligned&\Vert \sum_{\xi'\in Q}(\tilde\phi,\xi')\check\xi(t)\Vert_{H^s}\\
=&[\sum_j|j|^{2s}| \sum_{\xi'\in Q}(\tilde\phi,\xi')(\check\xi(t),\delta_j)|^2]^{1/2}\\
=&[\sum_j|j|^{2s}| \sum_k\sum_{\xi'\in Q}\hat\phi(k)\xi'(k, 0)(\check\xi(t),\delta_j)|^2]^{1/2}\endaligned\tag 4.28$$
From the support of $\xi'$ (2.25), 
$$||j|-|k||\leq 2(\log T)^\sigma\, (\sigma>2).\tag 4.29$$ 
Since $$|j|>2J_0=8A(\log T)^{\sigma}\quad (A>1)\tag 4.30$$
from (4.3), (4.29, 4.30) imply
$$|j|/2<|k|<2|j|.\tag 4.31$$

We now make a dyadic decomposition of $\phi$.  Let $R$ be dyadic and 
$$R/2<|j|<2R.\tag 4.32$$
So $$R/4<|k|<4R.\tag 4.33$$
Let $$\phi_R=\sum_{R/4<|k|<4R}\hat\phi(k)e^{ikx}.\tag 4.34$$
We then have 
$$\aligned(4.28)&\leq[\sum_{R \text{ dyadic}}4^sR^{2s}\sum_{R/2<|j|<2R}| \sum_k\sum_{\xi'\in Q}\hat\phi_R(k)\xi'(k, 0)(\check\xi(t),\delta_j)|^2]^{1/2}\\
&\leq [\sum_{R \text{ dyadic}}4^sR^{2s}\Vert \sum_{\xi'\in Q}(\tilde\phi_R, \xi')\check\xi(t)\Vert_2^2]^{1/2}\endaligned\tag 4.35$$
Using (4.27) and since $\text {supp }\phi_R\subset\text{ supp }\phi\subseteq[-J/2, -2J_0]\cup[2J_0, J/2]$,
$$\aligned\Vert \sum_{\xi'\in Q}(\tilde\phi_R, \xi')\check\xi(t)\Vert_2&\leq \Vert S(t)\phi_R\Vert_2+
2e^{-\frac{1}{2}(\log T)^{(\frac{\sigma'-\delta}{\alpha})}}\Vert\phi_R\Vert_2\\
&\leq 2\Vert\phi_R\Vert_2\quad (\sigma'>2\alpha+\delta>2).\endaligned\tag 4.36$$
Using (4.36) in (4.35), we have 
$$\aligned &\Vert \sum_{\xi'\in Q}(\tilde\phi, \xi')\check\xi(t)\Vert_{H^s}\\
\leq &[\sum_{R \text{ dyadic}}4^sR^{2s}\cdot 4\Vert\phi_R\Vert_2]^{1/2}\\
\leq&C^s\Vert\phi\Vert_{H^s}.\endaligned\tag 4.37$$
Combining (4.37) with (4.27, 4.1), we obtain (4.9) with a slightly larger $C$. $\square$ 
\enddemo
\smallskip
\demo{Proof of the Theorem}
We use the decomposition in (4.4-4.6), which decomposes into low, intermediate and high frequencies.
(4.9) controls (4.5), the intermediate frequencies: $2J_0\leq |j|\leq J/2$, (4.7) controls (4.6), the high
frequencies: $|j|>J/2$. So the only work left is to control (4.4), the low frequencies:
$|j|\leq 2J_0$, which we do by iterating $S(0, T){\overset\text{def }\to=}S(T)$, $|T|$ times and each time making again the decomposition as in (4.4-4.6).

We have 
$$\align&\Vert\Pi_{J/4}S(0, t)\Pi_{2J_0}u_0\Vert_{H^s}\\
\leq &\Vert\Pi_{J/4}S(1, t)\Pi_{2J_0}S(0, 1)\Pi_{2J_0}u_0\Vert_{H^s}\tag 4.38\\
+&\Vert\Pi_{J/4}S(1, t)(\Pi_{J/2}-\Pi_{2J_0})S(0, 1)\Pi_{2J_0}u_0\Vert_{H^s}\tag 4.39\\
+&\Vert\Pi_{J/4}S(1, t)(I-\Pi_{J/2})S(0, 1)\Pi_{2J_0}u_0\Vert_{H^s},\tag 4.40\endalign$$
which is the analogue at $t=1$ of the decomposition in (4.4-4.6), with $S(0, 1)\Pi_{2J_0}u_0$
replacing $u_0$. So we have 
$$\aligned(4.39)&\leq C^s\Vert S(0, 1)\Pi_{2J_0}u_0\Vert_{H^s}\\
&\leq 2C^{2s}s!,\endaligned\tag 4.41$$
where we used $$\Vert S(0, 1)\Vert_{H^s\to H^s}\leq 2C^ss!\tag 4.42$$ from (3.3) and 
$$\aligned (4.40)&\leq \Vert S(0, 1)\Pi_{2J_0}u_0\Vert_{H^s}\cdot \frac{(C^{s}s!)^4}{J}(1+|T|^{3s+1})\\
&\leq 2C^ss!\endaligned\tag 4.43$$
(4.41, 4.42) are the analogues of (4.9, 4.7), which control (4.5, 4.6).

Using (4.41, 4.42), we have after one iteration:
$$\align&\Vert\Pi_{J/4}S(0, t)\Pi_{2J_0}u_0\Vert_{H^s}\\
\leq &\Vert\Pi_{J/4}S(1, t)\Pi_{2J_0}S(0, 1)\Pi_{2J_0}u_0\Vert_{H^s}\\
&+4C^{2s}s!.\tag 4.44\endalign$$
After $r$ iterations, the analogue of the bound on (4.39) is 
$$\aligned &C^s\Vert S(r, r-1)\Pi_{2J_0}S(r-1, r-2)\Pi_{2J_0}\cdots\Pi_{2J_0}u_0\Vert_{H^s}\\
&\leq C^s \Vert S(r, r-1)\Vert_{H^s\to H^s}\cdot (2J_0)^s\\
&\leq 2 C^{2s}s! (2J_0)^s;\endaligned\tag 4.45$$
while the analogue of the bound on (4.40) is 
$$2C^ss!(2J_0)^s.\tag 4.46$$
After $|T|$ iterations, we then have 
$$\aligned& \Vert\Pi_{J/4}S(0, t)\Pi_{2J_0}u_0\Vert_{H^s}\\
\leq&\Vert\Pi_{2J_0}S(T-1, T)\Pi_{2J_0}\cdots\Pi_{2J_0}
S(r-1, r)\Pi_{2J_0}\cdots\Pi_{2J_0}S(0, 1)\Pi_{2J_0}u_0\Vert_{H^s}\\
&\,+4C^{2s}s!(2J_0)^s|T|\\
\leq &|T|(sJ_0)^s\cdot C^s\endaligned\tag 4.47$$
with a larger $C$.

Using (4.47) in (4.4) and combining with (4.9, 4.7), we obtain 
$$\Vert\Pi_{J/4}S(0, T)u_0\Vert_{H^s}\leq |T|(sJ_0)^s C^s.\tag 4.48$$
Using (4.48) in (4.2), we have 
$$\Vert S(0, T)u_0\Vert_{H^s}\leq C^s|T|(sJ_0)^s\tag 4.49$$
for all $0<s\leq\log T$. 
Interpolating with the $L^2$ bound $\Vert u_0\Vert_{L^2}\leq 1$ yields
$$\aligned \Vert S(0, T)\Vert_{H^{s'}\to H^{s'}}&\leq |T|^{s'/s}(CsJ_0)^{s'}\\
&\leq C^{s'}(\log T)^{(\sigma+1)s'}\quad (\sigma>2)\endaligned\tag 4.50$$
with a larger $C$, for all $0<s'<s$, where we took $s=\log |T|$ and used (4.3). 

For a fixed $s>0$, for $|t|<e^s$, the a priori bound (3.3) gives 
$$\Vert S(0, t)\Vert _{H^s\to H^s}\leq C^s s! e^{s^2},\tag 4.51$$
for $|t|\geq e^s$, we use (4.50). This gives immediately
$$\Vert S(0, t)\Vert _{H^s\to H^s}\leq C_s(\log (|t|+2))^{(\sigma+1)s}$$ for all $s>0$.
Let $\varsigma=\sigma+1$, we obtain the theorem. $\square$
\enddemo

\enddocument
\end